\documentclass[reqno]{amsart}
\usepackage{amsmath,amsthm,amssymb}
\usepackage{latexsym}
\usepackage{eucal}

\newtheorem{theorem}{Theorem}[section]
\newtheorem{lemma}{Lemma}[section]

\theoremstyle{definition}
\newtheorem{definition}{Definition}[section]

\def\cal{\mathcal}
\let\Re=\undefined
\DeclareMathOperator{\Re}{Re}
\let\Im=\undefined
\DeclareMathOperator{\Im}{Im}

\begin{document}
\title[Absolutely continuous spectrum of  multidimensional \ldots]{Absolutely continuous spectrum of  multidimensional Schr\"odinger
operator}
\author{Sergey A. Denisov}
%\address{Mathematics, 253-37, Caltech, Pasadena, CA, 91125, USA,
%{\rm e-mail: denissov@caltech.edu}}
\maketitle

\begin{abstract}
We prove that 3-dimensional Schr\"odinger operator with slowly
decaying potential has an a.c. spectrum that fills $\mathbb{R^+}$.
Asymptotics of Green's functions is obtained as well.
\end{abstract} \vspace{1cm}

 Consider the Schr\"odinger operator
\begin{equation}
H=-\Delta+V, x\in \mathbb{R}^d \label{hamiltonian}
\end{equation}
We are interested in finding the support of an a.c. spectrum of
$H$ for the slowly decaying potential $V$. The following
conjecture is due to B. Simon \cite{rarry}

{\bf Conjecture.} If $V(x)$ is such that
\begin{equation}
\int\limits_{\mathbb{R}^d} \frac{V^2(x)}{1+|x|^{d-1}}\, dx<\infty
\label{barry}
\end{equation}
then $\sigma_{ac}(-\Delta+V)=\mathbb{R}^+$.

It was proved for the one-dimensional case by Deift and Killip
\cite{dk} (see also \cite{rowan, mnv, denis7}). For some Dirac
operators, this conjecture was shown to be true for $d=1$ by M.
Krein \cite{Krein} and for $d=3$ by the author \cite{denis8}. For
the Schr\"odinger operator, certain multidimensional  results were
obtained recently in \cite{lns, lns1, kla}. The spatial
asymptotics of the Green function is a classical subject
\cite{Ag1,Ag2}. In the current paper we deal with $d=3$ and prove
the preservation of an a.c. spectrum under more restrictive
conditions on the potentials rather than (\ref{barry}). Methods of
the paper can be generalized to other $d$ and perhaps to the
discrete case too. We take $d=3$ for simplicity only.
 The structure of the paper is as follows. In the introduction, we
 obtain different results
 that serve as a motivation to the main theorem of the paper. In the second section,
 we prove the spatial asymptotics of the Green kernel. Then, in
 the third part, the main result on the preservation of the a.c.
 spectrum is obtained. The last section contains different
 applications.

Let us introduce some notations we will be using later. We denote
the integral kernel of $R_z=(H-z)^{-1}$ by $G_z(x,y)$. Recall that
for the Green's kernel of the free Laplacian, we have
\begin{equation}
G_z^0(x,y,z)=\frac{\exp (ik|x-y|) }{4\pi |x-y|}, z=k^2, k\in
\mathbb{C}^+
\end{equation}
The symbol $\Sigma$ stands for the unit sphere in $\mathbb{R}^3$.
The inner product of two vectors $\xi,\zeta$ in $\mathbb{R}^3$ is
denoted by $\langle \xi,\zeta \rangle$.
\newpage

{\bf Acknowledgement.} We are grateful to B. Simon for the useful
discussions.

\section{ Introduction}

We will start with the following simple result.

{\bf Proposition.} {\it Consider the $C^1(\mathbb{R}^3)$
vector-field $Q(x)$ such that
\begin{equation}
\int\limits_{\mathbb{R}^3} \frac{|Q(x)|^2}{1+|x|^2}\, dx<\infty
\end{equation}
and  $div\; Q(x)\in L^\infty(\mathbb{R}^3)$. Let $V=\gamma \cdot
div\; Q+|Q|^2, |\gamma|\leq 1$ and $H=-\Delta+V$. Then, $
\sigma_{ac}(H)=\mathbb{R}^+ $. }

 The proof of this fact follows
immediately from the arguments given in \cite{lns}. Integrating by
parts in the quadratic form for $H$, one can easily see that
$H\geq 0$ for any $|\gamma|\leq 1$. Then, since
\[
\left| \;\int\limits_{\mathbb{R}^3}
\frac{V(x)}{1+|x|^2}\,dx\right|<\infty
\]
the first trace-inequality (6.8) from \cite{lns}  yields
$\sigma_{ac}(H)=\mathbb{R}^+$. Since $H\geq 0$, one does not have
to worry about the negative eigenvalues and an analysis in
\cite{lns} is now easy.

Analogous argument works for any $d$ including $d=1$. But in the
one-dimensional case one then can argue that $Q^2$ is a relative
trace-class perturbation and the Rosenblum-Kato theorem \cite{rs}
would yield $\sigma_{ac}(-d^2/dx^2+Q')=\mathbb{R}^+$ for $Q$ being
any $L^2(\mathbb{R})$ function (see \cite{denis9}). In the
meantime, this argument does not work in the multidimensional
case. Even assuming  $|Q(x)|<C/(1+|x|)^{0.5+\varepsilon},
\varepsilon>0$, one has $|Q(x)|^2$-- short-range only and the
trace-class argument does not work (see \cite{yafaev}, p.22,
Problem 2.12). Still, we will consider this case and apply
different technique to show the preservation of the a.c. spectrum.
But first we want to discuss the following problem. In the
one-dimensional case, the positivity of the operators
$H_{\pm}=-d^2/dx^2 \pm Q'+Q^2$ on $\mathbb{R}$ follows from the
following factorization identity

\[
D=\left[
\begin{array}{cc}
0 & d/dx+Q\\
-d/dx+Q & 0
\end{array}
\right],
D^2=
\left[
\begin{array}{cc}
H_+ & 0\\
0 & H_-
\end{array}
\right]
\]
Operator $D$ corresponds to a certain Krein system \cite{Krein},
which simply makes the one-dimensional scattering theory a branch
of the approximation theory, in particular, the theory of
orthogonal polynomials. In $d>1$ case, we don't know analogous
result. Still one can come up with the following substitute.
Consider the following operators
\[
L=\left[
\begin{array}{cccc}
0 & -\partial_{x_1} & -\partial_{x_2} & -\partial_{x_3}\\
\partial_{x_1} & 0 & -\partial_{x_3} & \partial_{x_2}\\
\partial_{x_2} & \partial_{x_3} & 0 & -\partial_{x_1}\\
\partial_{x_3} & -\partial_{x_2} & \partial_{x_1} & 0\\
\end{array}
\right], M_v= \left[
\begin{array}{cccc}
0 & -v_1 & -v_2 & -v_3 \\
 v_1 & 0 & v_3 & -v_2 \\
 v_2 & -v_3 & 0 & v_1 \\
 v_3 & v_2 & -v_1 & 0
\end{array}
\right]
\]
acting in, say, $[S(\mathbb{R}^3)]^4$ and $v(x)$ is some real,
smooth vector-field. We introduce
\[
\cal{D}=\left[
\begin{array}{cc}
0 & L+M_v \\
L-M_v & 0
\end{array}
\right]
\]
The straightforward calculations show that $\cal{D}^2$ has $(1,1)$
component in the block representation equal to
$H=-\Delta+|v|^2+div \;v$. Thus, $H\geq 0$. If $v=\nabla \nu$,
than $V=\Delta \nu +|\nabla\nu|^2$ and the other elements in the
first raw/column of $\cal{D}^2$ are all zero, i.e. $\cal{D}^2$ is
a direct sum of the scalar Schr\"odinger operator and another
operator whose characterization is rather complicated. Infact, the
operator $\cal{D}$ has a
 very nice structure. Consider the unitary matrix $Y$
\[
Y=\left[
\begin{array}{cc}
0 & U\\
U & 0
\end{array}
\right], U=\frac{1}{\sqrt 2}\left[
\begin{array}{cccc}
1 & 0 & 0 & -i\\
0 & -i & 1 & 0\\
0 & -i & -1 & 0\\
1 & 0 & 0 & i
\end{array}
\right]
\]
Then, $Y\cal{D}Y^{-1}$ has the ``Dirac operator" form that allows
one to apply the methods of the paper \cite{denis8}. In
particular, it show that the Schr\"odinger operator with the
potential $V=\Delta \nu+|\nabla \nu|^2$ has the Green kernel with
certain spatial asymptotics as long as $|\nabla
\nu(x)|<C/(1+|x|)^{0.5+\varepsilon}, \varepsilon>0$. In this
paper, we choose a direct method to study Schr\"odinger operators
with potentials being the divergence of a slowly-decaying
vector-field. Notice here that this type of potentials was studied
earlier in the papers \cite{mv1, mv2}.

\section{ Asymptotics of the Green's function}

Let $0<\delta<C$ be fixed.  We begin with the following auxiliary
results
\begin{lemma}
Assume that $1<\rho<2|x|/3$. Then

\begin{equation}
\int\limits_{|y|=\rho }e^{-\delta (|x-y|+|y|)}d\tau _{y} <C\delta
^{-1}\rho e^{-\delta|x|} \label{lemma-ll1-1}
\end{equation}
\begin{equation}
\int\limits_{|y|=\rho }e^{-\delta(|x-y|+|y|)} \zeta (x,y)d\tau
_{y} <C\delta^{-1.5}  \rho^{0.5} e^{-\delta|x|}
\label{lemma-ll1-2}
\end{equation}
\begin{equation}
\int\limits_{|y|=\rho }e^{-\delta(|x-y|+|y|)} \zeta^2 (x,y)d\tau
_{y} <C\delta^{-2} e^{-\delta |x|}
\label{lemma-ll1-3}
\end{equation}
  \label{ll1}
\end{lemma}

Proof. Without loss of generality, assume $x=(0,0,|x|)$.
Introducing the spherical coordinates $y_{1}=\rho \cos \theta \cos
\varphi ,y_{2}=\rho
\cos \theta \sin \varphi ,y_{3}=\rho \sin \theta ,$ we get%
\begin{eqnarray*}
&&\rho ^{2}\int\limits_{-\pi }^{\pi }d\varphi \int\limits_{-\pi
/2}^{\pi /2}d\theta \cos \theta \exp \left( -\delta \left[ \rho
+\sqrt{|x|^{2}+\rho ^{2}-2|x|\rho \sin \theta }\right]\right) \\
&<&C\rho ^{2}e^{-|x|}\int\limits_{-\pi /2}^{\pi /2}d\theta \cos
\theta \exp \left[ -c\delta|x|\rho (|x|-\rho )^{-1}(1-\sin \theta
)\right] <Ce^{-\delta|x|}\frac{|x|-\rho }{\delta |x|}\rho
\end{eqnarray*}%
The estimate (\ref{lemma-ll1-1}) is now straightforward. To prove
(\ref{lemma-ll1-2}) and (\ref{lemma-ll1-3}), it suffices to notice
that $\zeta(x,y)\sim \sin \zeta(x,y)$ for small
$\zeta(x,y)$.$\blacksquare$

Let $|x|>1$ and $\Upsilon=\{y: |y|>2|x|/3, |x-y|>2|x|/3\}$.
\begin{lemma} The following estimate holds
\begin{equation}
\int\limits_\Upsilon \exp\left(-\delta \left[|x-y|+|y|\right]
\right) dy \leq C\delta^{-3} \exp(-\gamma \delta |x|)
\end{equation}
with $\gamma>1$. \label{lemma-ll2}
\end{lemma}
The proof repeats the one of lemma 3.4 from \cite{denis8} and is
elementary.

Consider the following class of functions in $\mathbb{R}^3$.
\begin{definition}
We say that $\psi(x)\in Cl(k), k\in \mathbb{C}^+$ if
$\psi(x)=\exp(ik|x|)\left[\psi_1(x)+\psi_2(x)\right]$, where
$\psi_{1(2)}$-- measurable, and

\[
\begin{array}{cc}
\displaystyle |\psi_1(x)|<\frac{1}{|x|^{1.5}+1}, &  {\rm Group }\; (A)\\
\displaystyle |\psi_2(x)|<\frac{1}{|x|+1}, |\nabla \psi_2 (x)|<
\frac{1}{|x|^{1.5}+1}, & {\rm Group }\;  (B)
\end{array}
\]
\end{definition}
Introduce an operator

\begin{equation}
B(k) f(x)=\int\limits_{\mathbb{R}^3} \frac{ \exp(ik|x-y|)}{4\pi
|x-y|} div[Q(y)] f(y)dy \label{operator-1}
\end{equation}

\begin{theorem}
Assume that $Q(x)$ is a vector-field such that
\[
|Q(x)|<\frac{m(Q)}{1+|x|^{0.5+\varepsilon}},\,
|div[Q(x)]|<\frac{m(Q)}{1+|x|^{0.5+\varepsilon}},\, \varepsilon>0
\]
Then, for
\[
|\Re k|<a,\, 0<\Im k<b,\, m(Q)<C(\varepsilon,a,b)[\Im k]^{3}
\]
$B(k)$ acts within the class $Cl(k)$. \label{th-asympt}
\end{theorem}

Proof. We will always assume that $k: |\Re k|<a, 0<\Im k<b$.
Denote $\Im k=\delta$.\\
 Group (A). Consider $\psi(x)=\exp(ikx)\psi_1(x)$. Let us show
that $B(k)\psi$ belongs in group (B). For $x: |x|<1$, all
estimates are easy. Assume that $|x|>1$. By lemma \ref{lemma-ll2},
we have

\[
\left|\, \exp(-ik|x|) \int\limits_{y\in \Upsilon}
\frac{\exp(ik|x-y|)}{|x-y|} div[Q(y)]
\exp(ik|y|)\psi_1(y)dy\right| <C(a,b)\delta^{-3} m(Q)(1+|x|)^{-3}
\]
Then, using (\ref{lemma-ll1-1}), we get

\[
\left| \,\exp(-ik|x|) \int\limits_{|y|<2|x|/3}
\frac{\exp(ik|x-y|)}{|x-y|} div[Q(y)]
\exp(ik|y|)\psi_1(y)dy\right|
\]
\[
 <\frac{C(a,b)m(Q)}{\delta (1+|x|)}
\int\limits_0^{2|x|/3} (1+\rho)^{-1-\varepsilon} d\rho\leq
\frac{C(a,b)m(Q)}{\delta (1+|x|)} \quad {\rm sharp\, !}
\]
Similarly, making change of variables $y-x=t$, we have
\[
\left|\, \exp(-ik|x|) \int\limits_{|y-x|<2|x|/3}
\frac{\exp(ik|x-y|)}{|x-y|} div[Q(y)]
\exp(ik|y|)\psi_1(y)dy\right|
\]
\[
=\left|\, \exp(-ik|x|) \int\limits_{|t|<2|x|/3}
\frac{\exp(ik|t|)}{|t|} div[Q(t+x)]
\exp(ik|t+x|)\psi_1(t+x)dt\right|
\]
\[
 <\frac{C(a,b)m(Q)}{\delta (1+|x|)^{2+\varepsilon}}
\int\limits_0^{2|x|/3}  d\rho\leq \frac{C(a,b)m(Q)}{\delta
(1+|x|)^{1+\varepsilon}}
\]
Let us show that the gradient of $\exp(-ik|x|)B(k)\psi$ is small.
Differentiating $|x-y|^{-1}$ gives a strong decay

\[
\left| \,\exp(-ik|x|) \int\limits_{|y|<2|x|/3} \frac{x-y}{|x-y|^3}
\exp(ik|x-y|) \cdot div[Q(y)] \cdot \exp(ik|y|) \cdot
\psi_1(y)dy\right|
\]
\[
<\frac{C(a,b)m(Q)}{\delta (1+|x|^2)} \int\limits_0^{2|x|/3}
(1+\rho)^{1+\varepsilon} d\rho< \frac{C(a,b,
\varepsilon)m(Q)}{\delta (1+|x|^2)};
\]
\[
\left| \, \exp(-ik|x|) \int\limits_{|y-x|<2|x|/3}
\frac{\exp(ik|x-y|)}{|x-y|^2} div[Q(y)]
\exp(ik|y|)\psi_1(y)dy\right|
\]
\[
<\frac{C(a,b)m(Q)}{\delta (1+|x|^{2+\varepsilon})}
\int\limits_0^{2|x|/3} (1+\rho)^{-1} d\rho<
\frac{C(a,b,\varepsilon)m(Q)}{\delta (1+|x|^2)}
\]
The integral over $\Upsilon$ is less than
$C(a,b)m(Q)\delta^{-3}(1+|x|^{-4})$.

The gradient of the exponent yields
\[
\left| \; \int\limits_{|y|<2|x|/3}
\frac{\exp(ik|x-y|-ik|x|)}{|x-y|} \left( ik \frac{x-y}{|x-y|}-ik
\frac{x}{|x|}\right) div[Q(y)] \exp(ik|y|)\psi_1(y)dy\right|
\]
\begin{equation}
< C(a,b) \left| \; \int\limits_{|y|<2|x|/3}
\frac{\exp(ik|x-y|-ik|x|)}{|x-y|} \zeta(x-y,x) div[Q(y)]
\exp(ik|y|)\psi_1(y)dy\right| \label{ineq1}
\end{equation}
 For
$\zeta(\xi_1,\xi_2)<\pi/2$, $\zeta(\xi_1,\xi_2)\sim \sin
\zeta(\xi_1,\xi_2)$. Due to $\sin$-- theorem, $\sin\zeta(x-y,
x)=\sin \zeta(x,y)|x-y|^{-1}|y|$. By (\ref{lemma-ll1-2}),
(\ref{ineq1}) is less than
\[
\frac{C(a,b)m(Q)}{\delta^{1.5}(1+|x|^2)}\int\limits_0^{2|x|/3}
\rho^{-0.5-\varepsilon} d\rho\leq \frac{C(a,b)m(Q)}{\delta^{1.5}
(1+|x|)^{1.5+\varepsilon}}
\]
Then,
\[
\left| \; \int\limits_{|y-x|<2|x|/3}
\frac{\exp(ik|x-y|-ik|x|)}{|x-y|} \zeta(x-y,x) div[Q(y)]
\exp(ik|y|)\psi_1(y)dy\right|
\]
\[
<\frac{C(a,b)m(Q)}{\delta^{1.5}(1+|x|^{2+\varepsilon})}
\int\limits_0^{2|x|/3} \rho^{-0.5}
d\rho<\frac{C(a,b)m(Q)}{\delta^{1.5}(1+|x|^{1.5+\varepsilon})}
\]
The integral over $\Upsilon$ is smaller than
$C(a,b)m(Q)\delta^{-3}(1+|x|^{-3})$. Thus, we showed that
$B(k)\psi(x)$ falls into group (B) assuming the corresponding
estimate on $m(Q)$.

Group (B). Consider $\psi(x)=\exp(ik|x|)\psi_2(x)$. Integrating by
parts, we have $\exp(-ik|x|) B(k)\psi(x)=I_1+I_2+I_3$. $I_1$ is
the term with the derivative falling onto $|x-y|^{-1}$, etc.

\[
|I_1|<\int\limits_{\mathbb{R}^3} |Q(y)| \frac{\exp\left[ \delta
(|x|-|x-y|-|y|)\right]}{|x-y|^2} |\psi_2(y)|dy
\]
Simple estimates show that
\[
|I_1(x)|<\frac{C(a,b,\varepsilon)m(Q)}{\delta^{3} (1+|x|^{1.5})}
\]
Therefore, $\exp(ik|x|)I_1(x)$ belongs to group (A) for $m(Q)$
sufficiently small. Let us show that $\exp(ik|x|)I_{2(3)}(x)$ are
in the group (B). For $I_2$,

\[
\int\limits_{|y|<2|x|/3} |Q(y)| \frac{\exp\left[ \delta \left(
|x|-|y|-|x-y|\right) \right]}{|x-y|} |\nabla \psi_2(y)|dy
\]
\[
<\frac{C(a,b)m(Q)}{\delta (1+|x|)}\int\limits_0^{2|x|/3}
(1+\rho)^{-1-\varepsilon}
d\rho<\frac{C(a,b,\varepsilon)m(Q)}{\delta (1+|x|)};
\]
\[
\int\limits_{|y-x|<2|x|/3} |Q(y)| \frac{\exp\left[ \delta \left(
|x|-|y|-|x-y|\right) \right]}{|x-y|} |\nabla \psi_2(y)|dy
\]
\[
<\frac{C(a,b)m(Q)}{\delta
(1+|x|^{2+\varepsilon})}\int\limits_0^{2|x|/3}
d\rho<\frac{C(a,b)m(Q)}{\delta (1+|x|^{1+\varepsilon})}
\]
The integral over $\Upsilon$ is smaller than
$C(a,b)m(Q)\delta^{-3}(1+|x|^{-3})$. Let us estimate $|\nabla
I_2(x)|$. It has two terms. The one containing derivative of
$|x-y|^{-1}$ is easy to deal with. It provides, again, the
stronger decay at infinity. The other term with the derivative of
the exponent can be bounded as follows

\[
\int\limits_{|y|<2|x|/3} \frac{|Q(y)|}{|x-y|} \cdot |\nabla
\psi_2(y)| \cdot \Big| \nabla_x \exp \left[ik\left( |x-y|+|y|-|x|
\right)\right]\Big| \;dy
\]
\[
<C(a,b)\int\limits_{|y|<2|x|/3} \frac{|Q(y)|}{|x-y|} |\nabla
\psi_2(y)|\; \zeta(x-y, x) \exp \left[-\delta\left(
|x-y|+|y|-|x|\right) \right]  dy
\]
\[
<\{{\rm by\; the\; \sin
-theorem}\}<\frac{C(a,b)m(Q)}{\delta^{1.5}(1+|x|^2)}
\int\limits_0^{2|x|/3} (1+\rho)^{-0.5-\varepsilon} d\rho<
\frac{C(a,b)m(Q)}{\delta^{1.5}(1+|x|^{1.5+\varepsilon})}
\]

\[
\int\limits_{|y-x|<2|x|/3} \frac{|Q(y)|}{|x-y|} |\nabla \psi_2(y)|
\Big| \nabla_x \exp \left[ik\left( |x-y|+|y|-|x|
\right)\right]\Big| dy
\]
\[
<\frac{C(a,b)m(Q)}{\delta^{1.5}(1+|x|^{2+\varepsilon})}
\int\limits_0^{2|x|/3} (1+\rho)^{-0.5} d\rho<
\frac{C(a,b)m(Q)}{\delta^{1.5}(1+|x|^{1.5+\varepsilon})}
\]
The integral over $\Upsilon$ is smaller than
$C(a,b)m(Q)\delta^{-3}(1+|x|^{-3})$. Thus, $\exp(ik|x|) I_2(x)$ is
in the group (B).

For $I_3$,

\[
\int\limits_{|y|<2|x|/3}  |Q(y)| \cdot |\psi_2(y)| \cdot
\zeta(y,x-y)
 \frac{\exp\left[ \delta \left( |x|-|y|-|x-y|\right)
\right]}{|x-y|} dy
\]
\[
<\{{\rm by\; the\; \sin -
theorem}\}<\frac{C(a,b)m(Q)}{\delta^{1.5}
(1+|x|)}\int\limits_0^{2|x|/3} (1+\rho)^{-1-\varepsilon}
d\rho<\frac{C(a,b,\varepsilon)m(Q)}{\delta^{1.5} (1+|x|)};
\]

\[
\int\limits_{|y-x|<2|x|/3}  |Q(y)| \cdot |\psi_2(y)| \cdot
\zeta(y,x-y)
 \frac{\exp\left[ \delta \left( |x|-|y|-|x-y|\right)
\right]}{|x-y|} dy
\]
\[
<\frac{C(a,b)m(Q)}{\delta^{1.5}
(1+|x|^{1.5+\varepsilon})}\int\limits_0^{2|x|/3} (1+\rho)^{-0.5}
d\rho<\frac{C(a,b)m(Q)}{\delta^{1.5} (1+|x|^{1+\varepsilon})};
\]
The integral over $\Upsilon$ is smaller than
$C(a,b)m(Q)\delta^{-3}(1+|x|^{-2.5})$.

Take the gradient of $I_3(x)$. Differentiation of
\[
\frac{1}{|x-y|}\left( \frac{y}{|y|}-\frac{x-y}{|x-y|} \right)
\]
in $x$ gives the term which can be estimated in the standard way
by \mbox{$C(a,b,\varepsilon) \delta^{-3} (1+|x|)^{-1.5}$}. Taking
the derivative of the exponent, we have
\[
\int\limits_{|y|<2|x|/3}  |Q(y)| \cdot |\psi_2(y)| \cdot
\zeta(y,x-y)\cdot \zeta(x,x-y)
 \frac{\exp\left[ \delta \left( |x|-|y|-|x-y|\right)
\right]}{|x-y|}\; dy
\]
\[
<\{{\rm by\; the\; \sin- theorem}\}<\frac{C(a,b)m(Q)}{\delta^{2}
(1+|x|^2)}\int\limits_0^{2|x|/3} (1+\rho)^{-0.5-\varepsilon}
d\rho<\frac{C(a,b)m(Q)}{\delta^{2} (1+|x|)^{1.5+\varepsilon}};
\]

\[
\int\limits_{|y-x|<2|x|/3}  |Q(y)| \cdot |\psi_2(y)| \cdot
\zeta(y,x-y)\cdot \zeta(x,x-y)
 \frac{\exp\left[ \delta \left( |x|-|y|-|x-y|\right)
\right]}{|x-y|} dy
\]
\[
<\frac{C(a,b)m(Q)}{\delta^{2}
(1+|x|)^{1.5+\varepsilon}}\int\limits_0^{2|x|/3} (1+\rho)^{-1}
d\rho<\frac{C(a,b,\varepsilon)m(Q)}{\delta^{2} (1+|x|)^{1.5}};
\]
The integral over $\Upsilon$ is smaller than
$C(a,b)m(Q)\delta^{-3}(1+|x|)^{-2.5}$. Thus, for small $m(Q)$,
$\exp(ik|x|) I_3$ falls in group (B) and the proof is
finished.$\blacksquare$

Remark. Clearly, the estimates for the integrals over $\Upsilon$
can be improved.

The following two theorems provide an asymptotics of the Green
function. Fix any $\varepsilon
>0$.
\begin{theorem}
Assume that $Q(x)$ is a vector-field such that
\[
|Q(x)|<\frac{m(Q)}{1+|x|^{0.5+\varepsilon}},\,
|div[Q(x)]|<\frac{m(Q)}{1+|x|^{0.5+\varepsilon}}
\]
Take $z=k^2, k=\tau+i\delta, 0<a_1<\tau<a_2, 0<\delta<b$. Let
$V=div\; Q$ in (\ref{hamiltonian}).

If
\begin{equation}
\delta^3>C(a_1,a_2,b) m(Q) \label{smallness}
\end{equation}
then
\begin{equation}
|G_z(x,y)-G_z^0(x,y)|<\frac{C(a_1,a_2,b)
m(Q)}{\delta^3-C(a_1,a_2,b) m(Q)} \frac{\exp (-\delta |x|)}{|x|}
\label{est-1}
\end{equation}
uniformly for $|y|<1, |x|>1$. \label{theorem-31}
\end{theorem}

Proof. Fix $\varepsilon>0$ and then $a_1, a_2, b$.  Consider the
second resolvent identity for $H$:
\[
(H-z)^{-1}=(H_0-z)^{-1}-(H_0-z)^{-1}V(H-z)^{-1}
\]
Therefore,
\begin{equation}
G_z(x,y)=G_z^0(x,y)-B(k)G_z(\cdot,y) \label{second-res}
\end{equation}
with $B(k)$ introduced in (\ref{operator-1}). Consider $y$ as a
parameter, $|y|<1$. Iterate (\ref{second-res}). It is an easy
exercise to show that $B(k) G_z^0(\cdot,y)\in Cl(k)$ for small
$m(Q)$. Therefore, (\ref{est-1}) follows directly from the theorem
\ref{th-asympt} by summing up the geometric series. $\blacksquare$

\begin{theorem}
Assume, again, that $Q(x)$ is a vector-field such that
\[
|Q(x)|<\frac{C}{1+|x|^{0.5+\varepsilon}},\, |div \;
Q(x)|<\frac{C}{1+|x|^{0.5+\varepsilon}}
\]
Take $z=k^2, k=\tau+i\delta, 0<a_1<\tau<a_2, 0<\delta<b$. Let
$V=div\,Q$ in (\ref{hamiltonian}). Consider any $f(x)$ with the
support inside the unit ball and $\|f\|_2<1$. Let
$u(x,k)=(H-z)^{-1}f$. The following estimate is true
\begin{equation}
\limsup_{|x|\to\infty} \Big\{ |x| \exp (\delta |x|)\cdot |u(x,k)|
\Big\} \leq A(\delta)  \label{est-2}
\end{equation}
and
\begin{equation}
A(\delta)\leq
\exp\left[C(a_1,a_2,b)\delta^{-\gamma(\varepsilon)}\right]
\end{equation}
with $\gamma(\varepsilon)>0$. \label{th-asa}
\end{theorem}

Proof. Fix any $\delta>0$.  We can control the Green function only
if $m(Q)$ in the theorem \ref{theorem-31} is relatively small. The
idea now is to cut out the big ball of radius $R(\delta)$ to
guarantee that $m(Q)$ satisfies assumptions of the theorem
\ref{theorem-31}, but relative to different $\varepsilon$, say
$\varepsilon/2$. Take $R>0$, it will be assigned with the precise
value later. Consider the radially-symmetric function $\chi_R(x)$:

\[
\chi_R(x)=\left\{
\begin{array}{cr}
1, & {\rm if}\; |x|<R;\\
0, & {\rm if}\; |x|>R+1
\end{array}
\right.
\]
We can always assume that $|\nabla\chi_R|<C$, where $C$ is
independent of $R$. Write $V=V_1+V_2$ where
\[
V_{1(2)}=div\; Q_{1(2)}, Q_1=\chi_R Q, Q_2=(1-\chi_R) Q
\]
For $Q_2$,
\begin{equation}
|Q_2(x)|<\frac{CR^{-\varepsilon/2}}{1+|x|^{0.5+\varepsilon/2}},
|div\,
Q_2(x)|<\frac{CR^{-\varepsilon/2}}{1+|x|^{0.5+\varepsilon/2}}
\end{equation}
Consider $H_2=-\Delta+ V_2$. We have $u(x,k)=(H_2-z)^{-1}[ f-V_1 u
]$. Denote the resolvent kernel of $H_2$ by $L_z(x,y)$. Then,

\begin{equation}
|u(x,k)|<C\int\limits_{|y|<R+1} |L_z(x,y)|\cdot |u(y,k)|dy <
C\|u(\cdot,k)\|_2 \Big[ \int\limits_{|y|<R+1}
|L_z(x,y)|^2dy\Big]^{0.5} \label{lll-1}
\end{equation}
Clearly, $\|u(\cdot,k)\|_2 <C\delta^{-1}$.

Fix any $y, |y|<R+1$. The function $L_z(x,y)=K_{z,y}(x-y,0)$ where
$K_{z,y}(s,t)$ is the resolvent kernel of the Schr\"odinger
operator with the shifted potential $V_y(s)=V_2(s+y)$. Notice that
\[
|V_y(s)|<\frac{CR^{-\varepsilon/2}}{(1+|s|)^{0.5+\varepsilon/2}}
\]
where $C$ is independent of $y, |y|<R+1$. Let
$\varepsilon_1=\varepsilon/2$. Now, fix $R$ such that (see
(\ref{smallness}))
\[
C(a_1,a_2,b) R^{-\varepsilon_1}<\delta^3
\]
We can take $R=[\delta^3/(2C(a_1,a_2,b) ]^{-1/\varepsilon_1}$.
Then, the theorem \ref{theorem-31} is applicable and we have (see
(\ref{est-1}))
\[
 |x|\exp(\delta |x|)\cdot |K_{z,y}(x,0)| <C(a_1,a_2,b),
|y|<R+1
\]

Thus
\[
\limsup\limits_{|x|\to\infty} \Big\{|x| \exp (\delta |x|) \cdot
|L_z(x,y)|\Big\}< C(a_1,a_2,b)\exp[C\delta R], |y|<R+1
\]
and, by (\ref{lll-1}),
\begin{equation}
\limsup_{|x|\to\infty} \Big\{ |x| \exp (\delta |x|)\cdot |u(x,k)|
\Big\} \leq C\delta^{-1} R^{1.5} \sup_{|y|<R+1}
\limsup_{|x|\to\infty} \Big\{ |x| \exp (\delta |x|) \cdot
|L_z(x,y)|\Big\}
\end{equation}
Consequently,
\[
\limsup_{|x|\to\infty} \Big\{ |x| \exp (\delta |x|)\cdot |u(x,k)|
\Big\}<\exp [C\delta^{-\gamma(\varepsilon)}]
\]
with $\gamma(\varepsilon)>>1$. $\blacksquare$

Remark 1. Modifying slightly the proof of the theorem
\ref{th-asympt}, one should be able to drop the condition $|div\;
Q(x)|<C(|x|+1)^{-0.5-\varepsilon}$. The mere boundedness of $V$
could be enough. The other results of the paper then should also
follow.

Remark 2. Perhaps, one can work directly with the equation
$-\Delta u+Vu-zu=f$ to improve an estimate (\ref{est-2}). We
expect at most polynomial growth for $A(\delta)$ as $\delta \to
0$. It is a reasonable guess that $\limsup_{|x|\to\infty}
|x|\exp(\delta |x|) \cdot \|G_z(|x|\cdot
\theta,0)\|_{L^2(\theta\in\Sigma)}$ is finite under the condition
(\ref{barry}).

\section{ Absolutely continuous spectrum}

We will start with an easy but fundamental factorization identity.
Consider $H$ with compactly supported bounded potential $V(x)$.
Take any $f(x)\in L^\infty(\mathbb{R}^3)$ with a compact support.
Let $\Pi$ be a rectangle in $\mathbb{C}^+$: $k=\tau+i\delta,
0<a_1<\tau <a_2, 0<\delta<b$, and $u(x,k)=(H-z)^{-1}f, z=k^2$. We
have
\begin{equation}
\begin{array}{cc}
\displaystyle u(x,k)=&\displaystyle \frac{\exp (ikr)}{r}\left(
A(k,\theta)+\bar{o}(1)\right), \\
\displaystyle \frac{\partial u(x,k)}{\partial r}= &ik
\displaystyle \frac{\exp
(ikr)}{r}\left( A(k,\theta)+\bar{o}(1)\right),\\
\displaystyle & r=|x|, \theta=\displaystyle \frac{x}{|x|},
|x|\to\infty
 \end{array}
 \quad {\rm (Sommerfeld's\; radiation\; conditions)}
\end{equation}

The amplitude $A(k,\theta)$ has the following properties:

\begin{itemize}
\item  $A(k, \theta)$ is an analytic in $k\in \Pi$
vector-function.

\item The absorption principle holds, i.e. $A(k,\theta)$ is
continuous on $\overline\Pi$.

\item For the boundary value of the resolvent, we have
(\cite{yafaev}, p.40--42)
\\ \phantom \quad \mbox{$\Im (R^+_{k^2}f,f)=k\|A(k,\theta)\|_{L^2(\Sigma)}^2,
k>0$}. Therefore,
\begin{equation}
\sigma'_f(E)=k\pi^{-1}\|A(k,\theta)\|_{L^2(\Sigma)}^2, E=k^2
\label{factor}
\end{equation}
where $\sigma_f(E)$ is the spectral measure of $f$.
\end{itemize}

The main result of this section is the following multidimensional
version of the corollary on p.181, \cite{denis9}.

\begin{theorem}
Let $Q(x)$ be a vector-field in $\mathbb{R}^3$ and
\[
|Q(x)|<\frac{C}{1+|x|^{0.5+\varepsilon}},\, |div\;
Q(x)|<\frac{C}{1+|x|^{0.5+\varepsilon}}, \varepsilon>0
\]
Then, $H=-\Delta+div\; Q$ has an a.c. spectrum that fills
$\mathbb{R}^+$. \label{ac-spectrum}
\end{theorem}

Proof. Let us fix any interval $I=[a_1,a_2]\subset R^+$ and $b>0$.
We will show that $I\subset \sigma_{ac}(H)$. Following
\cite{rowan}, consider an isosceles triangle $T$ in $\Pi$ with the
base equal to $I$ and the adjacent angles both equal to
$\pi/\gamma_1$, $\gamma_1>\gamma(\varepsilon)$ with
$\gamma(\varepsilon)$ from (\ref{est-2}).

Take $f(x)$ as any nonzero $L^\infty(\mathbb{R}^3)$ function
supported on the unit ball. Then
\[
A_0(k,\theta)=\lim_{|x|\to\infty}
|x|\exp(-ik|x|)\int\limits_{\mathbb{R}^3} \frac{\exp
[ik|x-y|]}{4\pi |x-y|}f(y)dy
\]
\[
=(4\pi)^{-1} \int\limits_{|y|<1} \exp[-ik\langle \theta,y\rangle
]f(y)dy
\]
For the fixed $\theta'$, $A_0(k,\theta')$ is entire in $k$.
Therefore, we can find a point $k_0=\tau_0+i\delta_0$ inside the
triangle $T$ such that $A_0(k_0,\theta')\neq 0$. Since
$A_0(k_0,\theta)$ is continuous in $\theta$,
\begin{equation}
\|A_0(k_0,\theta)\|_{L^2(\theta\in \Sigma)}>0 \label{choice-k}
\end{equation}
Fix this $k_0$ for the rest of the proof. Consider $R>0$ and the
function $\chi_R(x)$ introduced in the proof of the theorem
\ref{th-asa}. Let, again, $Q_1=\chi_R Q, Q_2=(1-\chi_R) Q$ and
$V_{1(2)}=div\; Q_{1(2)}$. Notice that $V_1$ is compactly
supported. Therefore, by the trace-class argument,
$\sigma_{ac}(-\Delta+V)=\sigma_{ac}(-\Delta+V_2)$. Thus, we can
restrict our attention to $H_2=-\Delta+V_2$ only. For $V_2$, we
have $|Q_2(x)|<CR^{-\varepsilon/2}/(1+|x|)^{0.5+\varepsilon/2}$
and $|div\;
Q_2(x)|<CR^{-\varepsilon/2}/(1+|x|)^{0.5+\varepsilon/2}$. Take $R$
big enough to have the following estimate
\begin{equation}
\limsup_{|x|\to \infty} \Big\|\; |x|\exp (-ik_0 |x|)\cdot
(H_2-k_0^2)^{-1} f\Big\|_{L^2(\Sigma)}>0 \label{amplit-1}
\end{equation}
We can always do that due to theorem \ref{theorem-31} and
(\ref{choice-k}). Fix this $R$ and the corresponding $V_2$.

 Now, take any $\rho>R+1$ and consider $Q^{(\rho)}=\chi_\rho Q_2, V^{(\rho)}=div\; Q^{(\rho)}$.
Since $V^{(\rho)}$ is compactly supported, an amplitude $A_\rho
(k,\theta)$ of $f(x)$ is well defined. Consider the following
function $\nu_\rho (k)=\ln \|A_\rho (k,\theta)\|_{L^2(\Sigma)}$.
It is subharmonic in $T$. Let $\omega(k_0,s), s\in \partial T$
denote the value at $k_0$ of the Poisson kernel associated to $T$.
One can easily show that

\begin{equation}
0\leq \omega(k_0,s)<C|s-s_{1(2)}|^{\gamma_1-1},\;  s\in
\partial T \label{bound}
\end{equation}
where $s_{1(2)}$ are endpoints of $I$. By subharmonicity,
\begin{equation}
\int\limits_{s\in \partial T} \nu_\rho (s) \omega(k_0,s) d|s|\geq
\nu_\rho (k_0)
\end{equation}
If we denote the edges of $\partial T$ by $I_{1(2)}$, i.e.
$\partial T=I\cup I_1\cup I_2$, then

\begin{equation}
\int\limits_{s\in I} \omega(k_0,s) \ln \|A_\rho
(s,\theta)\|_{L^2(\Sigma)} ds\geq \nu_\rho(k_0)-\int\limits_{s\in
I_{1}\cup I_{2}} \omega(k_0,s) \ln^+ \|A_\rho
(s,\theta)\|_{L^2(\Sigma)}ds
\end{equation}
Notice now that (\ref{amplit-1}), theorem \ref{th-asa}, and
(\ref{bound}) yield

\begin{equation}
\int\limits_{s\in I} \omega(k_0,s) \ln \|A_\rho
(s,\theta)\|_{L^2(\Sigma)} ds>C
\end{equation}
with the constant $C$ independent of $\rho$. Recall the
factorization identity (\ref{factor}) to have
\begin{equation}
\int\limits_{k\in I} \omega(k_0,k) \ln \sigma_{(f,\rho)}'(k^2)
dk>C
\end{equation}
where $\sigma_{(f,\rho)}$ is the spectral measure of $f$ with
respect to $-\Delta+V^{(\rho)}$. It is an easy exercise to show
that $(-\Delta+V^{(\rho)}-z)^{-1}\to (-\Delta+V_2-z)^{-1}$ in the
strong sense as $\rho\to\infty$. Therefore, $d\sigma_{(f,\rho)}\to
d\sigma_f$ in the weak-($\ast$) sense. The usual argument with the
semicontinuity of the entropy \cite{ks} then implies
\begin{equation}
\int\limits_{{J^2}} \ln \sigma_f'(E) dE>-\infty
\end{equation}
for any subinterval $J\subset I$. Since $I$ was an arbitrary
interval in $\mathbb{R}^+$, we have $\sigma_{ac}(H)=\mathbb{R^+}$.
$\blacksquare$

The case of radially-symmetric potential shows that a very rich
singular spectrum is allowed under the conditions of the theorem
\cite{denis9, kls}.

Remark. Consider the short-range potential $V(x):$
\begin{equation}
|V(x)|<C/(1+|x|^{1+\varepsilon}),\, \varepsilon>0 \label{short-r}
\end{equation}
Then one can easily show that the analogs of theorem
\ref{th-asympt}, \ref{theorem-31}, \ref{th-asa} hold. That allows
us to show that $\sigma_{ac} (-\Delta+V_1+V_2)=\mathbb{R}^+$ where
$V_1$ is from the theorem \ref{ac-spectrum} and $V_2$ satisfies
(\ref{short-r}).

\section{ Applications}

In this section, we consider some concrete examples of the
potentials.

 Example 1. Consider

\[
V(x)=\frac{\sin x_1}{(1+x_1^2+x_2^2+x_3^2)^\gamma}, \gamma>1/4
\]
One can write
\[
V(x)=-\frac{\partial}{\partial x_1} \frac{\cos
x_1}{(1+x_1^2+x_2^2+x_3^2)^\gamma}+V_2(x)
\]
with $V_2(x)$-- short-range. By remark after the theorem,
$\sigma_{ac}(-\Delta+V)=\mathbb{R}^+$.

 Then, by theorem \ref{theorem-31}, the asymptotics of the Green function
\mbox{$G_{-k^2}(x,0)$} is similar to the asymptotics of
$G_{-k^2}^0(x,0)$. Here $k$ is assumed to be sufficiently large.
In the meantime, take a point $x$ on the hyperplane $x_1=const$
very far from the origin. One can see that the Agmon distance
(i.e. solution to the corresponding eikonal equation) from $x$ to
the origin will have the WKB-type correction to the usual linear
growth. Thus, it is not the solution to the eikonal equation that
governs the phase of the Green function.\vspace{1cm}

 Talking about the asymptotics of Green's function, one can suggest
the following approach which works sometimes. Let us try to find
$u(x,k)=\exp(-k|x|+\mu(x))/|x|$ that solves $-\Delta u+Vu+k^2u=0$
for $|x|>1$. We also assume that $k\in \mathbb{R}^+$ and $k>>1$.
The equation for $\mu$ now reads as follows

\begin{equation}
\Delta \mu +|\nabla \mu|^2-2k\frac{\partial \mu}{\partial
r}=V+\frac{2}{r}\cdot \frac{\partial \mu}{\partial r}, \; r=|x|
\label{HJ}
\end{equation}
It is an eikonal equation with viscosity, modified by the radial
derivative term. Making the following substitution
$\mu=r\exp(r)\psi$, one ends up with a simple equation, which
yields
\begin{equation}
\mu=-GV+G[|\nabla\mu|^2] \label{aster}
\end{equation}
where the operator $G$ is defined as follows
\[ Gf(x)=
|x|\exp(k|x|)\int\limits_{\mathbb{R}^3}
\frac{\exp[-k(|x-y|+|y|)]}{4\pi |x-y|\cdot |y|} f(y)dy
\]
If $V$ is such that the gradient of $GV$ is decaying fast, then,
one can hope to iterate (\ref{aster}). Then the leading term in
the asymptotics of the phase $\mu$ would be $GV$. Unfortunately,
this idea does not work unless we assume the strong decay of the
derivatives of $V$. As the first example suggests, neither the
nonlinear term, nor the viscosity can be discarded in (\ref{HJ}).
We do not know the right WKB correction to the asymptotics of the
Green kernel for potentials $|V(x)|<C/(1+|x|)^{0.5+\varepsilon}$.
That is a major problem that prohibits us from proving Simon's
conjecture in its original form. An advantage of the equation
(\ref{aster}) is that it contains the potential in $GV$ form only.
This function $GV$ is an integral and provides a lot of averaging
for $V$. That averaging might be useful for studying the random
Schr\"odinger operators.\vspace{1cm}

 Example 2. Consider a smooth vector-field $Q(x)$ supported on the
 unit ball. Take
 \[
 V(x)=\sum\limits_{j\in \mathbb{Z}^+} a_j v(x-x_j)
 \]
 where $V(x)=div \; Q(x)$, points $x_j$ are scattered in
 $\mathbb{R}^3$ such that $|x_k-x_l|>2, k\neq l$, and $a_j\to 0$
 such that $|V(x)|<C/(1+|x|^{0.5+\varepsilon})$. Then the theorem \ref{ac-spectrum} is
 applicable and $\sigma_{ac}(-\Delta+V)=\mathbb{R}^+$.\vspace{1cm}

  So far, we were able to deal with potentials that can be written
 as the divergence of a slowly-decaying field. But on the formal
 level, any function can be written as a divergence of some
 vector-field, for instance

 \[
 V(x)=\Delta \Delta^{-1} V=-div\; \nabla_x \int\limits_{\mathbb{R}^3}
 \frac{V(y)}{4\pi |x-y|} dy=div\; \int\limits_{\mathbb{R}^3}
 \frac{x-y}{4\pi |x-y|^3} V(y)dy
 \]
One can easily show that for continuous $V$ with compact support,
this identity holds true. In the general situation, given $V(x)$,
one might consider the vector-field
\begin{equation}
Q(x)=\int\limits_{\mathbb{R}^3}
 \frac{x-y}{4\pi |x-y|^3} V(y)dy \label{potencial}
\end{equation}
try to show that $Q(x)$ is well defined, satisfies the bound
$|Q(x)|<C/(1+|x|^{0.5+\varepsilon})$, and $V=div\; Q$. Then, as
long as $V(x)$ itself is slowly-decaying, i.e.
$|V(x)|<C/(1+|x|^{0.5+\varepsilon})$, the theorem
\ref{ac-spectrum} would yield $\sigma_{ac}
(-\Delta+V)=\mathbb{R}^+$.\vspace{1cm}

Example 3 (The Anderson model with slow decay). Consider the
following model. Take a smooth function $\phi(x)$ with the support
inside the unit ball. Like in the second example, consider

\[
 V_0(x)=\sum\limits_{j\in \mathbb{Z}^+} a_j \phi(x-x_j)
 \]
 where the points $x_j$ are scattered in
 $\mathbb{R}^3$ such that $|x_k-x_l|>2, k\neq l$, and $a_j\to 0$
 in a way that $|V_0(x)|<C/(1+|x|^{0.5+\varepsilon})$. Let us now
 ``randomize" $V_0$ as follows

 \begin{equation}
 V(x)=\sum\limits_{j\in \mathbb{Z}^+} a_j \xi_j \phi(x-x_j)
 \label{defin}
 \end{equation}
 where $\xi_j$ are real-valued, bounded, independent random variables with \mbox{$\mathbb{E}
 \Big[ \xi_j^{2k+1}\Big] =0$}, $k\in \mathbb{Z}^+$. Clearly, any even distribution satisfies the last condition.

 \begin{theorem}
For $V$ given by (\ref{defin}), we have
$\sigma_{ac}(-\Delta+V)=\mathbb{R^+}$ almost surely.
 \end{theorem}

Proof.  Fix any $x_0$. By Kolmogorov's one series theorem, the
integral in (\ref{potencial}) converges almost surely, i.e.

\begin{equation}
\int\limits_{|y|<R}
 \frac{x_0-y}{4\pi |x_0-y|^3} V(y)dy \to Q(x_0), R\to\infty
 \label{justif-1}
\end{equation}
for $\omega\in \Omega$, $\mathbb{P}[\Omega]=1$.

For $x$ inside a fixed compact $K$, we use the Lagrange theorem to
have
\begin{equation}
 \left|
 \frac{x_0-y}{|x_0-y|^3}-\frac{x-y}{|x-y|^3}\right|<\frac{C(K)}{1+|y|^3},
 |y|>>1\label{justif-2}
\end{equation}
Therefore,
\[
Q(x)=\lim\limits_{R\to\infty} \int\limits_{|y|<R}
 \frac{x-y}{4\pi |x-y|^3}\, V(y)dy
\]
\[
 =
Q(x_0)+\lim\limits_{R\to\infty} \frac{1}{4\pi} \int\limits_{|y|<R}
 \left[ \frac{x-y}{ |x-y|^3}-\frac{x_0-y}{ |x_0-y|^3}\right] V(y)dy
 \]
exists for any $\omega\in \Omega$ and is continuous in $x\in
\mathbb{R}^3$. One also has $div\; Q(x)=V(x)$. Let us now show
that $|Q(x)|<C/(1+|x|)^{0.5+\varepsilon_1}$ with probability one
for some $\varepsilon_1>0$. We can write
\[
Q(x)=Q_1(x)+Q_2(x)=\int\limits_{|x-y|<1} \frac{x-y}{4\pi |x-y|^3}
V(y)dy+\int\limits_{|x-y|>1} \frac{x-y}{4\pi |x-y|^3} V(y)dy
\]
Clearly, for $Q_1(x)$ : $|Q_1(x)|<C/(1+|x|)^{0.5+\varepsilon}$. As
about $Q_2(x)$, we don't have singularity under the integral
anymore and one can easily show that
\begin{equation}
|D Q_2(x)|<C\ln(1+|x|)/(1+|x|)^{0.5+\varepsilon} \label{unifi}
\end{equation}
where $D$ means the differential of any component of $Q_2$. We
introduce
\[
S_j(x)=\int\limits_{|x-y|>1} \frac{x-y}{4\pi |x-y|^3}
\,\phi(y-x_j)dy
\]
Then,
\begin{equation}
\mathbb{E} \Big[ |Q_2(x)|^2 \Big] \leq \sum\limits_{j\in
\mathbb{Z}^+} a_j^2 \cdot \mathbb{E} [\xi_j^2] \cdot |S_j(x)|^2
\label{dispersiya}
\end{equation}
\[
< C \sum\limits_{j\in \mathbb{Z}^+} \frac{a_j^2}{1+|x-x_j|^4}< C
\int\limits_{\mathbb{R}^3}\frac{dy}{(1+|y|^{1+2\varepsilon})(1+|x-y|^4)}<\frac{C
}{1+|x|^{1+2\varepsilon}}
\]
Now, consider the following sum
\[
\sum\limits_{k\in \mathbb{Z}^3} |k|^\gamma |Q_2(k)|^{2p}
\]
where $|k|^2=k_1^2+k_2^2+k_3^2$. Let us prove it converges almost
surely for the suitable choice of $\gamma>0$ and $p\in
\mathbb{N}$. We calculate the expectation
\[
\mathbb{E} \left[ \sum\limits_{k\in \mathbb{Z}^3} |k|^\gamma
|Q_2(k)|^{2p}\right]=\sum\limits_{k\in \mathbb{Z}^3} |k|^\gamma
\cdot \mathbb{E} \Big[|Q_2(k)|^{2p}\Big]
\]
Consider
\[
 \mathbb{E} \Big[\langle Q_2(k),Q_2(k)\rangle^p\Big]
\]
\begin{equation}
\leq \sum\limits_{j_1,\ldots,j_{p},m_1,\ldots,m_p}
a_{j_1}a_{m_1}\ldots a_{j_{p}}a_{m_p} \mathbb{E} \Big[
\xi_{j_1}\xi_{m_1}\ldots\xi_{j_{p}}\xi_{m_p}\Big]\langle
S_{j_1}(k),S_{m_1}(k)\rangle \ldots \langle
S_{j_{p}}(k),S_{m_p}(k)\rangle  \label{erist}
\end{equation}
Since all odd moments of $\xi_j$ are zero, $\mathbb{E} \Big[
\xi_{j_1}\xi_{m_1}\ldots\xi_{j_{p}}\xi_{m_p}\Big]$ is nonzero iff
the indices $j_1,\ldots,j_{p},m_1,\ldots, m_p$ coincide pairwise.
Therefore,
\[
\mathbb{E} \Big[\langle Q_2(k),Q_2(k)\rangle^p\Big]
<C(p)\sum\limits_{l_1,\ldots,l_{p}} a_{l_1}^2\ldots a_{l_{p}}^2
|S_{l_1}(k)|^2\ldots |S_{l_{p}}(k)|^2
\]
where $C(p)$ is a combinatorial factor. Thus, just like in
(\ref{dispersiya}),
\[
\mathbb{E} \Big[|Q_2(k)|^{2p}\Big]<C(p)\left[\sum\limits_{l}
a_{l}^2 |S_{l}(k)|^2\right]^p<C(p)/(1+|k|)^{p(1+2\varepsilon)}
\]
So,
\[
\mathbb{E} \left[ \sum\limits_{k\in \mathbb{Z}^3} |k|^\gamma
|Q_2(k)|^{2p}\right]<\infty
\]
as long as $\gamma=p(1+2\varepsilon)-3-\delta$, $\delta>0$. So,
the sequence $|k|^\gamma |Q_2(k)|^{2p}\in
\ell^1(\mathbb{Z}^3)\subset \ell^{\infty}(\mathbb{Z}^3)$ almost
surely. That means
\[
|Q_2(k)|<C(1+|k|)^{-\frac{\gamma}{2p}}=
\frac{C}{1+|k|^{0.5+\varepsilon-\frac{3+\delta}{2p}}}
\]
Taking $p$ big enough, we see that
$|Q_2(k)|<C(1+|k|)^{-0.5-\varepsilon_1}$ almost surely. Constant
$C$, of course, is random. But since $Q_2$ satisfies an estimate
(\ref{unifi}), we have

\begin{equation}
|Q_2(x)|<C/(1+|x|)^{0.5+\varepsilon_2},
0<\varepsilon_2<\varepsilon_1
\end{equation}
for all $x\in \mathbb{R}^3$ with probability one. Thus, the
theorem \ref{ac-spectrum} is applicable.$\blacksquare$

Remark. The assumption that the odd moments of $\xi_j$ are zeroes
can probably be dropped. In this case, one might have nonzero
contribution from factors like $\mathbb{E} \Big[ \xi_j^3\Big]$,
etc. in the sum (\ref{erist}). Perhaps, the corresponding terms
can be estimated as well. We do not pursue it here. The
representation of the potential $V$ as a divergence of
slowly-decaying vector-field is a multidimensional phenomena. In
dimension one, the argument does not work. Notice that we not only
found the support of an a.c. spectrum, but also proved an
asymptotics of the Green function for the spectral parameter in
the resolvent set. In the discrete setting, Bourgain \cite{bour1,
bour2} obtains stronger results for the Anderson model with slow
decay. See also the following paper \cite{rods}.

\end{document}